\documentclass[11pt]{amsart}
\usepackage[english]{babel}
\usepackage{amsmath}
\usepackage{amsthm}
\usepackage{amssymb}
\usepackage{amscd}
\usepackage{tabularx}

\def\sist{\left\{\begin{array}{l}} \def\sistt{\end{array}\right.}
\def\bb{\mathbb}
\def\la{\langle}
\def\ra{\rangle}

\newcommand{\pp}[1]{\mathbb{P}^{#1}}
\def\p3{\pp{3}}
\def\pi2{\pp{2}}
\def\p5{\pp{5}}

\def\pm{\pp{m}}
\def\pN{\pp{N}}
\def\Gr{\mathbb{G}} \def\Ta{\mathbb{T}}
\def\phi{\varphi}

\def\Mu{\mathrm{M}}
\begin{document}
\title[Varieties with degenerate Gauss image]{Some remarks
on varieties with degenerate
Gauss image}

\author[E.\ Mezzetti,
   O.\ Tommasi] {E.\ Mezzetti, O.\ Tommasi}
\address{Dipartimento di Scienze Matematiche,
Universit\`a di Trieste, Via Valerio
12/1, 34127 Trieste, Italy}
\email{mezzette@univ.trieste.it}

\address{Department of Mathematics, University of Nijmegen,
Toernooiveld, 6525 ED Nijmegen,
The Netherlands}
\email{tommasi@math.kun.nl}

\thanks{Work supported by funds of the University of Trieste (fondi
60\%) and  MIUR, project \lq\lq Geometria sulle variet\`a
algebriche''. The first Author acknowledges the hospitality of
Newton Institute of Cambridge, during the preparation of this
paper. }


\newtheorem{thm}{Theorem}[section]
\newtheorem{prop}[thm]{Proposition}
\newtheorem{lem}[thm]{Lemma}
\newtheorem{cor}[thm]{Corollary}

\theoremstyle{definition}

\newtheorem{defn}[thm]{Definition}
\newtheorem{notation}[thm]{Notation}
\newtheorem{example}[thm]{Example}
\newtheorem{conj}[thm]{Conjecture}
\newtheorem{prob}[thm]{Problem}
  \newtheorem{ex*}{Example}
\theoremstyle{remark}

\newtheorem{rem}[thm]{Remark}


\begin{abstract}

We consider projective varieties with degenerate Gauss image
whose focal hypersurfaces
are non-reduced schemes. Examples of this situation are provided by
the secant varieties of
Severi and Scorza varieties. The Severi varieties are moreover
characterized by a
uniqueness property.
\
\end{abstract}

\maketitle

\section{Introduction}\label{intro}

A classical theorem on surfaces
states that, if the tangent plane to
a ruled projective surface $S$
remains fixed along a general line of
the ruling, then $S$ is a developable surface, i.e., a cone or the
tangent developable of a curve.

In the case of higher dimensional varieties, this generalizes to the
problem of giving a
structure theorem for projective
varieties with degenerate Gauss image.
More precisely, let $X\subset\pN$ be
a   projective variety of dimension $n$. The Gauss map of $X$ is the
rational map
$\gamma$ from $X$ to the
Grassmannian $\Gr(n,N)$ of $n$-dimensional subspaces of $\pN$,
associating to a smooth point $x$
of $X$ the embedded tangent space $\Ta_xX$ to
$X$ at $x$. The Gauss image $\gamma(X)$ of $X$ is,
by definition, the closure of
$\gamma(X_{sm})$, where $X_{sm}$ is the smooth locus of $X$.

Clearly $\dim\gamma(X)\leq n$, and for \lq\lq general''
varieties equality
holds. Several terms have been used recently to denote the varieties such
that $\dim\gamma(X) < \dim X$: they are called
varieties with degenerate Gauss image by Landsberg (\cite {L99}),  developable
by Piontkowski (\cite{P}), tangentially degenerate
by Akivis - Goldberg (\cite{AG}). We will follow  Landsberg's convention.

Several general facts are well known for these varieties
(see for
instance the classical article \cite{S}). 
First of all, general fibres of $\gamma$ are linear. This
means that,
if the dimension of the Gauss image is $r<n$,
then the general fibre is a linear space of dimension $n-r$ 
  along which the tangent space is constant.
Secondly, if $X$
is not linear, then it is singular and its singular locus cuts a general
fibre along a codimension 1 subscheme.

The theory of foci for families of linear varieties well applies
   to the family of fibres of $\gamma$.
This allows us to consider and study
the focal locus on $X$ and on the fibre $\Lambda$ of $\gamma$,
which we will denote by $F$ and $F_\Lambda$ respectively. Both loci
have a natural scheme structure. 

Several classification results have recently been established for these
varieties (e.g. \cite{GH}, \cite{AGL}, \cite{MT}, \cite{P}, 
\cite{AG}) but a  general structure theorem is still missing.
In \cite{AG} a structure theorem is proved,
but under  rather strong assumptions on a general $F_{\Lambda}$,
for example that it is reduced or set-theoretically linear.
Moreover, Akivis and Goldberg state the problem of constructing
examples of projective varieties with degenerate Gauss image such
that $F_{\Lambda}$ is not reduced and $(F_{\Lambda})_{red}$ is not
linear (or proving that such varieties don't exist).

Piontkowski (\cite{P}) proved a uniqueness
theorem for developable
varieties $X$ with Gauss image of dimension 2, such that the focal locus
has codimension 2 in $X$. Those for which the focal conic is
integral are precisely the varieties of secant lines of
the Veronese surface $v_2(\pi2)$ or of cones over
it. Inspired by this example, we construct a series of examples where
$F_{\Lambda}$
is not reduced and $(F_{\Lambda})_{red}$ has
arbitrarily large degree. The examples are the secant varieties of Severi
varieties of dimension $>2$ and some natural
generalizations of them.

Zak's theorem on linear normality (\cite{Z83}) states that a smooth
non\-de\-ge\-ne\-rate 
$m$-dimen\-sio\-nal subvariety of $\pN$, with $N<3m/2+2$,
cannot be isomorphically
projected to $\pp{N-1}$. Zak also classified in \cite{Z86} the varieties
which are border cases in his theorem, and called them Severi varieties.
He proved that there are only 4 smooth examples: the Veronese surface $v_2
(\pi2)$ in
$\p5$, $\pi2\times\pi2$, $\Gr(1,5)$ and the non-classical variety $E$ of
dimension 16. For any of these  varieties $F$, its secant variety $X=SF$ is
a cubic hypersurface in $\pN$
with Gauss image isomorphic to $X$.

We prove in Example \ref{exSeveri} 
that on a general fibre of the
Gauss map for the secant variety of a Severi variety, the focal
scheme is a quadric of maximal rank with multiplicity  1, 2, 4, 8
respectively.

In Section \ref{sec:examples} 
we show that these examples can be
generalized in two ways, taking instead of the Severi varieties
  the following three series of varieties $Y$:
$v_2(\pm)$,
$\pm\times\pm$ and $\Gr(1, 2m+1)$, of dimensions $m$, $2m$, $4m$ respectively
(they are Scorza varieties according to Zak \cite{Z93}). If we consider their
higher secant variety
$S^{m-1}Y$ of $m$-secant $(m-1)$-planes, it turns out that they
are hypersurfaces of degree $m+1$, whose dual is $Y$ and whose singular
locus is $S^{m-2}Y$: it appears in each fibre with
multiplicity resp. 1, 2, 4. If we consider $SY$ instead,
we get on each fibre a  focal quadric with an arbitrarily high
multiplicity, more precisely,
$m-1$, $2(m-1)$, $4(m-1)$ respectively.

It is interesting to observe that the fourth Severi variety does
not generalize to give a complete class of examples
(\cite{L}).

We observe also that there is an upper bound on 
the codimension of the focal scheme and that there are
restrictions on the set of possible codimensions if we require
that $F$ is not reduced (see Theorem \ref{thm:codim2}). For \lq\lq
high'' codimension we prove in Theorem \ref{thm:unicity} some
uniqueness results, in the case in which the general focal
hypersurface is a quadric of maximal rank with multiplicity 2. We
extend in this way  the result found by Piontkowski in the case of
the secant variety of the Veronese surface.  Our results rely on
the fact, stated by Zak in \cite{Z4}, that if $F$ is a variety
whose general linear section 
is a Severi variety, then $F$ is a cone over a Severi variety.

The plan of the article is as follows. In \S \ref{sec:foci} we
recall 
the language of focal schemes for families of linear spaces, as
introduced in  \cite{MT}.
We also prove Proposition
\ref{thm:codim}, giving a lower bound, based on the codimension of
the focal scheme, for the multiplicity of the focal scheme. In \S
\ref{sec:gaps} we prove Theorems \ref{thm:codim2} and
\ref{thm:unicity}. Finally in \S \ref{sec:examples}, we give two
series of examples, showing that the focal hypersurfaces can have
arbitrarily high degree and multiplicity.

\subsection{Acknowledgment}
We would like to thank Prof. Fyodor Zak for some useful
suggestions and the referees for interesting remarks and
criticism.

\subsection{Notations and conventions}

In this paper a variety will be an integral closed subscheme of a
projective space over an
algebraically closed field $\mathbb{K}$, $char\ \mathbb{K} =0$.\\
If $\Lambda\subset \bb{P}^N$ is a projective linear subspace,
$\hat{\Lambda}\subset \mathbb K$ will denote the linear subspace
associated to $\Lambda$ such that $\Lambda=\bb{P}(\hat{\Lambda})$.
$T_xX$ will denote
the Zariski tangent space to the variety $X$ at its point $x$, while we will
denote by $\mathbb{T}_xX\subset\bb{P}^N$ the
embedded tangent space to $X$ at $x$.\\
We will always
use the same symbol to denote the points of a Grassmannian and the
corresponding linear subspaces.

\section{The focal scheme and its multiplicity}\label{sec:foci}

Let $X\subset\pN$ be a projective  variety of dimension $n<N$. Let
us assume that $X$ is covered by an $r$-dimensional family  of
linear spaces of dimension $k:=n-r$. Let $B$ be the subvariety of
the Grassmannian $\Gr(k,N)$ parametrizing that family.

Let us denote by $\mathcal{I}$ the incidence correspondence of
$B$, with the natural projections:
\begin{eqnarray}
\begin{array}
{ccccc} \;B&\xleftarrow[\mspace{50mu}]{p_1}&B \times\bb{P}^N&
\xrightarrow[\mspace{50mu}]{p_2}&\bb{P}^N\\[-6pt]
&&\cup &&\\[-6pt] B&\xleftarrow[\mspace{50mu}]{g}&\mathcal{I}&
\xrightarrow[\mspace{50mu}]{f}&\; X. \end{array}
\end{eqnarray}

We will associate to the family $B$  its focal subscheme
$\Phi\subset\mathcal{I}$ and its focal locus $F\subset X$. For
their definitions and for the proof of results cited in this
section, we refer to \cite{MT}.

We recall moreover that 
 the {\em characteristic map} of the family
$B$ is the map $\chi:=\beta\circ\alpha$,

$$\left(p_1^*(\mathcal{T}_{B})\right)\,_{|\mathcal{I}}
\xrightarrow{\ \alpha\ } \mathcal{T}_{\mathcal{I}} \xrightarrow{\
\beta\ } \mathcal{N}_{ \mathcal{I}|B \times \mathbb{P}^N},$$ where
$\alpha$ comes from the exact sequence expressing the tangent
sheaf 
to the product variety $B\times \pN$ as a product of tangent
sheaves, and $\beta$ from that defining the normal sheaf to
$\mathcal I$ inside $B\times \pN$.

For every smooth point $\Lambda$ of $B$, 
the restriction of $\chi$ to $g^{-1}(\Lambda)$ is called the
\emph{characteristic map of $B$ relative to $\Lambda$}.

Finally the {\em focal scheme on $\Lambda$}, denoted $F_\Lambda$,
is the scheme-theoretic intersection of the Cartier divisor
$\Phi\subset \mathcal{I}$ with
$\{\Lambda\}\times\Lambda\subset\mathcal I$.

\bigskip Let $\gamma:X\dasharrow \Gr(n,N)$ be the (rational) Gauss
map of $X$, regular on $X_{sm}$. It is a consequence of biduality
(\cite{K})
that its fibres are linear spaces, i.e. that a general
embedded tangent space is tangent to $X$ along a linear space. In
the following, we will apply the above contruction
to the family of the fibres of $\gamma$, assuming 
that they are of positive dimension.

\begin{defn}
The \emph{Gauss rank} of a variety $X$ is the number 
$r=\dim\gamma(X)$. If $r<n$, then $X$ is called a \emph{variety
with degenerate Gauss image}.
\end{defn}

Varieties with degenerate Gauss image are characterized by the
fact that the focal scheme 
 on a general $\Lambda\in B$ is a
hypersurface of degree $r$. Another important property is that the
focal locus of $B$ is always contained in the singular locus of
$X$.

Let us denote by $\bar g$ (respectively, $\bar f$) the restriction
of $g$ (resp., $f$) to $\Phi\subset\mathcal I$.
Note that the fibre $\bar{g}^{-1}(\Lambda)$ coincide with $F_\Lambda$,
which has in general dimension $n-r-1$, hence every irreducible component
of $\Phi$ has dimension $n-1$.

\begin{prop}\label{thm:codim} 
Let $X$ be a variety of dimension $n$ with Gauss image of
dimension $r<n$. Let $\bar\Phi$ be an irreducible component of the
focal scheme of the family $B$ considered as a reduced variety.
Let $\bar F\subset X$ be its scheme-theoretic image and let $c$ be
the codimension of $\bar F$ in $X$.
For a  general point $P=(\Lambda, x)\in\bar\Phi$, the Cartier
divisor $\Phi\subset\mathcal{I}$ has Samuel multiplicity
$\mu(P,\Phi)\geq c-1.$
\end{prop}

\begin{cor}\label{cor} For a general element $(\Lambda, x)\in\bar\Phi$,
either $F_\Lambda=\{\Lambda\}\times\Lambda$, or the Cartier
divisor $F_\Lambda\subset\{\Lambda\}\times\Lambda$ has Samuel
multiplicity $\mu(x,F_\Lambda)\geq c-1$.
\end{cor}
\proof Let $P=(\Lambda, x)\in \bar \Phi$ be a point that projects to $x$.
Considering the differentials of the maps $f$ and $\bar f$
at $P$, we get the diagram:
\begin{eqnarray}\begin{array} {ccc}\;T_P\mathcal{I}&\xrightarrow[
\mspace{50mu}]{d_Pf}&T_xX\\ [-6pt]\cup&&\cup\\[-6pt] T_P\bar\Phi&\xrightarrow[
\mspace{50mu}]{d_P\bar f}&\;T_x\bar F
\end{array}\end{eqnarray} We consider also the characteristic map
$\chi(\Lambda)$
of the family $B$ relative to $\Lambda$. At the
point $P$, it gives rise to a linear map $$\chi(\Lambda, x):
T_{\Lambda}B\rightarrow
\hat{\mathbb{T}}_xX/\hat{\Lambda}.$$

Let $K_P$ denote the common kernel of both $d_Pf$ and $\chi(\Lambda,
x)$. By generic
smoothness $d_P\bar f$ is surjective, hence $$\dim\ker d_P\bar
f=(n-1)-(n-c)=c-1.$$  But $\ker
d_P\bar f=K_P\cap T_P\bar\Phi$, where $T_P\bar\Phi$ has dimension $n-1$, so
either $K_P\subset T_P\bar\Phi$, the two kernels
coincide and $\dim K_P=c-1$, or they are different and $\dim K_P=c$.
Since the structure of scheme on $\Phi$ is given by the minors of
the characteristic map, this proves that $\mu(P,\Phi)$ is at least
$c-1$ in the former case, and at least $c$ in the latter.
    \qed
\vskip 12pt

\begin{ex*}
The multiplicity $\mu(x,F_\Lambda)$ can be strictly greater than $c$.
We see now an example of such a situation, in which the fibres of
the Gauss map have special properties of tangency to the focal
locus.

Let us
consider a birational map $\phi: F\dasharrow S$ between two
surfaces $F,S\subset\pp6$. Then
$$\begin{matrix}\psi:& F&\dasharrow& \Gr(3,6)\\
&x&\mapsto&\la \mathbb{T}_xF,\phi(x)\ra 
\end{matrix}$$
defines a rational map on $F$, 
 and $X$ is defined by the closure
of the variety swept by the 3-planes in $\pp{6}$ belonging to the
image of $\psi$. If the choice of $F,S$ and $\phi$ is general,
then $X$ is a variety of dimension 5 whose Gauss image has
dimension 4, and its focal locus
 is $F$, with codimension 3. If $x$ is a
general point of $F$,
then any line passing through $x$ and contained in $\la\mathbb{T}_xF,
\phi(x)\ra$ is a fibre of the Gauss map.

Let $P=(\Lambda, x)$ project to
a general point of $F$. By a direct computation, $x$ results to be 
the only focus on $\Lambda$, with multiplicity $4=c+1$. 

It is nonetheless true that $\dim K_P=2=c-1$:
the kernels of $d_P\bar f$ and $d_Pf$
indeed coincide and are generated by the directions of the curves in
$\mathcal{I}$ of the form $(\Lambda(t),x)$, where $\Lambda(t)$ varies in the star of lines of centre $x$ contained in $\psi(x)$.

The just constructed example belongs to the class of hyperbands of
\cite{AGL}. It can be generalized to a whole series of analogous
examples with larger dimensions and increasing difference between 
the multiplicity of the focal locus and the codimension $c$.

\end{ex*}

\begin{ex*}\emph{The Severi varieties.}\label{exSeveri}

Let $F$ denote one of the four Severi varieties (see the 
Introduction).
Let $X=SF$ be its secant variety, i.e. the closure of the 
union of lines joining two distinct points of $F$.
There are several known facts about $F$.
\begin{prop}[\cite{Z86}]
Let $F\subset\pN$ be a Severi variety. Then the following hold:
\begin{itemize}
\item[-] F has dimension $m\in\{2,4,8,16\}$;
\item[-] F is embedded in a projective space of dimension $N=3m/2+2$;
\item[-]
the secant variety $X$ of $F$ is a normal cubic hypersurface in $\pN$;
\item[-] $X$ is isomorphic to the dual variety of $F$;
\item[-] the singular locus of $X$ coincides with $F$.
\end{itemize}
\end{prop}
As $X$ is a hypersurface, its Gauss image is the same as its dual variety,
which is $F$. Then $X$ has Gauss rank $m$ and the fibres of the Gauss
map $\gamma$ of $X$
are linear subspaces of dimension
$m/2+1$. As shown in \cite{Z86}, the intersection with $F$ of the secant
lines of $F$ passing
through a general point $x\in X$ is a
quadric generating a space of dimension $m/2+1$. This space is
precisely the fibre
of the  Gauss map passing through $x$, and $B$
is just the family of such spaces. As the quadric depends on the
fibre $\Lambda$
to which $x$ belongs and not on $x$, we will
denote it by $Q_\Lambda$.

The degree of the focal scheme on $\Lambda$ must equal the Gauss rank,
so it is $m$, i.e. 2, 4, 8
or 16 respectively. But the focal locus of $B$
on $\Lambda$ has to be $Q_{\Lambda}$, of degree 2,
because it has to be contained in the singular locus of $X$. This
shows that $F_{\Lambda}$ coincides with $Q_{\Lambda}$ set-theoretically,
but, as a scheme,
  appears with multiplicity 1, 2,
4, 8 respectively in the four cases. This fact is confirmed by
Proposition \ref{thm:codim}, because in this case $c=m/2+1$ and
the multiplicity $\mu(x,\bar F)$ is equal to $m/2=c-1$.
\end{ex*}

\section{Gaps on the codimension and uniqueness of Severi examples}
\label{sec:gaps}
The results of the previous section allow us to state some bounds
on the possible codimension of the focal locus for the family of
the fibres of the Gauss map and some uniqueness results.

\begin{thm}\label{thm:codim2}
Let $X$ be a variety of dimension $n$ with degenerate Gauss image
$\gamma(X)$ of dimension $r$ with $r\geq 2$.  Let $c$ denote the
codimension of the focal locus $F$ in $X$. Then:
\begin{itemize}
\label{properties}
    \item[(i)]
$c\leq r+1$;
    \item[(ii)]  if
$c=r+1$, then $X$ is a cone of vertex a space
of dimension $(n-r-1)$
over a variety of dimension $r$;
    \item[(iii)]  if
    $c\leq r$ and the reduced focal hypersurfaces
    are not linear, then  $c\leq r/2+1$;
\item[(iv)] if $c=r/2+1$, then either the multiplicity of the
focal hypersurfaces in the fibres of $\gamma$
is at least $c$ and they are set-theoretically linear, 
or the multiplicity is $r/2$
    and the reduced focal hypersurfaces
are quadrics.
\end{itemize}
\end{thm}
\proof (i) Let $F_{\Lambda}$ be the focal scheme on a general
fibre $\Lambda$ of the Gauss map. Then $F_{\Lambda}$ is a
hypersurface of degree $r$, and from  Corollary \ref{cor},
it has multiplicity $\mu(x, F_\Lambda)\geq c-1$ at a
general point $x$. Hence $c-1\leq \mu(x, F_\Lambda)\leq r$, which implies
(i).

(ii) If $c=r+1$, then a general $F_{\Lambda}$ is a hyperplane in
$\Lambda$ with multiplicity $r$. On the other hand $\dim F=n-r+1=
\dim F_{\Lambda}$, so $F_{\Lambda}$ is a fixed $(n-r+1)$-plane and
$X$ is a cone over it.

(iii) and (iv) Assume $r/2+1\leq c\leq r$: then $\mu(x,F_\Lambda) \geq c-1
\geq r/2$.
So the degree of the reduced focal locus on $\Lambda$
is $r/\mu(x, F_\Lambda)\leq 2$.  If $\mu(x,F_\Lambda)\geq c$, then the hypothesis implies that
$\mu(x,F_\Lambda)=r$, so $(F_{{\Lambda})_{red}}$ is linear.
If $\mu(x,F_\Lambda)=c-1$, then our hypothesis implies $c= r/2+1$ and
$\mu(x,F_\Lambda)=r/2$. \qed
    \vskip 12pt

\begin{rem}
The description of varieties with Gauss rank $r=1$ is classical (see \cite{FP}
for a modern account), the description
of varieties with $r=2$
has been  accomplished recently by Piontkowski (\cite{P}). His uniqueness
theorem for the case in which the 
focal locus has codimension 2
can be extended as
follows. \end{rem}

\begin{thm}\label{thm:unicity}
Let $X$ be a variety of dimension $n$ with degenerate Gauss image $\gamma(X)$
    of dimension $r\geq 2$
and let $F$ be the focal locus of the family
of fibres of $\gamma$.  Assume that the
codimension of  $F$ in $X$ is $c=r/2+1$ and that the focal hypersurfaces
are quadrics counted with multiplicity $r/2$.
\begin{itemize}
\item[(i)] If $F$ is irreducible, then $X$ coincides with $SF$,
the secant variety of $F$,
and $n\geq3r/2+1$;
\item[(ii)] if moreover $n=3r/2+1$ and a general quadric $F_\Lambda$ is
smooth, then $F$ is a Severi variety;
\item[(iii)] if $n\geq 3r/2+1$, then the rank of the quadrics
$F_\Lambda$ is at most $r/2+1$, and, if equality generically holds, then
$F$ is (a cone over) a Severi variety
of dimension $r$. \end{itemize}
\end{thm}
\proof
(i) Since the focal hypersurfaces have set-theoretically degree 2, then the
lines which are contained in the fibres of $\gamma$ are all secant lines of
$F$. Hence $X\subset SF$. On the other hand, the family of lines
obtained in this way has 
dimension $r+\dim \Gr(1,n-r)$, which is equal to $2\dim F$. This proves
that $X=SF$.
A standard count of parameters shows that the family of quadrics passing
through a
general point of $F$ has dimension $r/2$ and that the intersection of two focal
quadrics
is always non-empty.

Let us denote by $B$ as usual the family of fibres of $\gamma$.
Let $\Sigma\subset B\times
B\times F$ be the set of triples $(\Lambda,\Mu,P)$ such that
$P\in F_\Lambda\cap F_\Mu$. Considering
the projections from $\Sigma$ to $B\times B$ and to
$F$, one gets easily that a general intersection $F_\Lambda\cap F_\Mu$ has
dimension
$n-3r/2-1$. This proves (i).

(ii) We can assume without loss of generality that $X$ is non-degenerate
(otherwise, we can restrict to the projective subspace of minimal
dimension in which $X$ is contained).
Suppose first that $X$ is a hypersurface.
We have that $F$ is a variety of dimension $r$ covered by a
$r/2$-dimensional family of non-singular quadrics of dimension $r/2$,
such
that a general pair of such quadrics meet at a point.
This property of the family of quadrics allows us to argue as in 
 \cite[Lemma 5]{Z86}, and deduce, as there,  that $F$ is smooth,
so it is a Severi variety by definition.

Suppose now that $X\subset\pp N$,
$N=3/2 r+2+k$, $k>0$. The image of $X$ in a projection with centre a
general $(k-1)$-plane of $\pp N$ is a 
hypersurface, which is the secant variety of a Severi variety $F'$
by the previous argument. Since the focal locus of $X$ project to
$F'$, we have that the general section of $F$ by a linear subspace
of dimension $3/2r+2$ is a Severi variety. By Corollary 1 of \cite{Z4}
Severi varieties are unextendable, hence $F$ must be a cone over
$F'$, which is impossible because in that case the quadrics
$F_\Lambda$ would be singular.

   (iii) Assume that the rank of $F_\Lambda$ is $\geq r/2 +1$.
Let $h:=(n-r/2-1)-(r+1)=n-3/2r-2$, $h\geq -1$. Cutting $F$ with
$L$, the intersection of $h+1$
general hyperplanes, we get $F':=F\cap L$, of dimension $r$,
containing a $r$-dimensional family of  quadrics of dimension $r/2$
and generically of maximal rank $r/2+1$.
Moreover the variety of secant lines
of $F'$ coincides with $SF\cap L$, so it has dimension $3r/2+1$.
By (ii) we obtain that $F'$ is a Severi variety.
Again since 
Severi varieties are unextendable, we can conclude that $F$ is
   a cone over $F'$ with vertex a
linear space of dimension $h$. In particular, the rank of $F_\Lambda$
is $r/2 +1$.
\qed

\begin{rem} Smooth varieties $F$
such that $X=S F$ satisfies the
assumptions of Theorem \ref{thm:unicity} are precisely the varieties studied by
Ohno in \cite{O}.

\end{rem}

\section{Examples of focal hypersurfaces of any degree or with arbitrarily
high multiplicity}\label{sec:examples}
In this section
we give two series of examples: the former show
that the focal
hypersurfaces in the fibres of the Gauss map can appear with
arbitrarily high multiplicity,
the latter have
focal hypersurfaces of arbitrary degree.
Note that in Example 3 the variety $X$ is not a hypersurface if
$m\geq 3$.

\begin{ex*} \emph{Secant varieties of Scorza varieties.}\label{exScorza1}
Let $F$ be one of the following varieties: $v_2(\pm)$, $\pm\times\pm$,
$\Gr(1,2m+1)$. They are contained in the projective spaces of dimension
${{m+2}\choose 2} -1$, $(m+1)^2-1$, ${{2m+2}\choose 2} -1$ respectively,
and have dimensions $m$, $2m$, $4m$. In all cases $F$ is defined by
suitable minors of
a matrix of variables, and precisely by $2\times 2$ minors of a symmetric
$(m+1)\times(m+1)$ matrix in the first case, 
by $2\times 2$ minors of a generic
$(m+1)\times(m+1)$ matrix
in the second case, by Pfaffians of $4\times 4$ minors of a
skew-symmetric matrix of
order $2m+2$ in the third one. These varieties are considered by Zak
in \cite{Z93},
and named by him Scorza varieties.

Let $X$ be the secant variety of $F$: it is defined by minors of the
same matrix as $F$
of order one more.
Then the fibres of the
Gauss map of $X$ have respectively dimension 2, 3, 5,
the focal hypersurfaces are quadrics
of dimensions 1, 2, 4, which appear with multiplicity $m-1$,
$2m-2$, $4m-4$ respectively.
\end{ex*}

\begin{ex*} \emph{Higher secant varieties of Scorza
varieties.}\label{exScorza2}
Let $X=S^{m-1}F$ be the  variety of
$m$-secant $(m-1)$-planes of $F$, 
where $F$ is one of the varieties appearing in the previous
example. $X$ is the maximal proper
secant variety of $F$ and it is the hypersurface defined
by the determinant (or the Pfaffian) of the matrix
considered in previous example.
The varieties $X$ and $F$ result to be mutually dual, 
so $r=m,2m,4m$ respectively.
The focal locus of the Gauss map is $S^{m-2}F$.
The focal hypersurfaces are the maximal secant varieties of
$v_2(\mathbb{P}^{m-1})$,
$\mathbb{P}^{m-1}\times\mathbb{P}^{m-1}$, $\Gr(1,2m-1)$
respectively.
They have degree $m$ and appear with multiplicity
respectively 1, 2, 4.
\end{ex*}

\begin{rem}
The class of Scorza varieties includes also
$\mathbb{P}^{m}\times\mathbb{P}^{m-1}$
and $\Gr(1,2m)$, for all $m\geq 3$.
Also their secant varieties have properties similar to above.
They are defined by minors of matrices of type $(m+1)\times (m+2)$
and of skew-symmetric square matrices
of type $2m+1$ respectively.

\end{rem}


\begin{thebibliography}{mmm}\addcontentsline{toc}{section}{\numberline{}
References}

\bibitem[AG]{AG}M. A. Akivis, V. V. Goldberg,
\emph{On the structure of submanifolds with degenerate Gauss maps},
Geom. Dedicata,
{\bf 86} (2001), 205-226

\bibitem[AGL]{AGL}M. A. Akivis, V . V. Goldberg, J.
M. Landsberg, \emph{Varieties with degenerate Gauss mapping},
preprint, 2000, math.AG/9908079

\bibitem[FP]{FP}G. Fischer, J. Piontkowski, \emph{Ruled surfaces},
Viehweg, 2001

\bibitem[GH]{GH}Ph. Griffiths, J. Harris, \emph{Algebraic Geometry and Local
Differential Geometry}, Ann. Sci. \'Ecole Norm. Sup. (4)
{\bf 12} (1979),  355-452

\bibitem[K]{K} S.L. Kleiman, 	\emph{Concerning the dual variety},  Progr. Math., 11, Birkh\" auser, Boston, 
    1981,  386--396,

\bibitem[L1]{L}J. M. Landsberg, \emph{On degenerate
secant and tangential varieties and local
differential geometry},
Duke Math. J. {\bf 85}, No.3 (1996), 605-634

\bibitem[L2]{L99}J. M. Landsberg, \emph{Algebraic Geometry and Projective
Differential Geometry},
Lecture Note Series, n. 45, Seoul National Univ., Seoul, Korea, 1999




\bibitem[MT]{MT}E. Mezzetti, O. Tommasi, \emph{On projective  varieties of
dimension $n+k$ covered by $k$-spaces}, Illinois J. Math., {\bf 46}, No.2 (2002), 443-465

\bibitem[O]{O}M. Ohno, \emph{On degenerate secant varieties whose Gauss maps
   have the largest images}, Pacific J. Math., {\bf 44}, No.1 (1999), 151-175

\bibitem[P]{P}J. Piontkowski, \emph{Developable varieties of Gauss
rank 2}, Internat. J. Math., {\bf 13}, No. 1 (2002), 93-110


\bibitem[S]{S}C. Segre, \emph{Preliminari di una teoria delle variet\`a luoghi di
spazi}, Preliminari di una teoria
delle variet\`a luoghi di spazi, Rend. Circ. Mate. Palermo (1)
30 (1910),  87-121, also in Opere Scelte, Vol. 2, ed. Cremonesi, Roma, 1958

\bibitem[Z1]{Z83}F.L. Zak, \emph{Projections of algebraic varieties}, Math.
USSR Sbornik, {\bf 44} (1983), 535-544

\bibitem[Z2]{Z86}F.L. Zak, \emph{Severi varieties},
Math. USSR Sbornik, {\bf 54} (1986), 113-127

\bibitem[Z3]{Z93}F.L. Zak, \emph{Tangents and secants of algebraic varieties},
Translations of Mathematical Monographs. 127. Providence, RI:
American Mathematical Society, 1993 

\bibitem[Z4]{Z4}F.L. Zak, \emph{Some properties of dual varieties and their applications
in projective geometry},  Lecture Notes in Math., 1479,
    Springer, Berlin, 1991, 273--280,
\end{thebibliography}
\end{document}